\documentclass{amsart}
\usepackage{amssymb,amscd}  

\newtheorem{theorem}{Theorem}[section]

\newtheorem{corollary}[theorem]{Corollary}

\newtheorem{proposition}[theorem]{Proposition}

\theoremstyle{definition}

\newtheorem{remark}[theorem]{Remark} 
\numberwithin{equation}{section}

\newcommand\C{\mathbf{C}} 
\newcommand\R{\mathbf{R}}
\newcommand\Q{\mathbf{Q}}
\newcommand\Z{\mathbf{Z}}

\newcommand\isomorphic{\cong}

\newcommand\Directsum{\bigoplus}

\DeclareMathOperator{\tr}{tr}

\newcommand\disjointunion{\sqcup}
\newcommand\intersect{\cap}

\newcommand\modforms{\mathcal{M}}
\newcommand\modring{{\modforms_*}} 

\newcommand\Half{\mathcal{H}} 
\newcommand\wtk{\kappa} 
\newcommand\levelN{N} 
\newcommand\modformsk{\modforms_\wtk} 
\newcommand\ring{\mathcal{R}}
\newcommand\ringo[1]{\ring_{\Gamma_0(#1)}}
\newcommand\ringi[1]{\ring_{\Gamma_1(#1)}}
\newcommand\ringoN{\ringo{\levelN}} 
\newcommand\GoN{\Gamma_0(\levelN)} 
\newcommand\GiN{\Gamma_1(\levelN)}
\newcommand\GpmN{\Gamma_{\pm}(\levelN)} 
\newcommand\TN{T_\levelN} 
\newcommand\troN{\tr^{\GoN}_{\Gamma(1)}} 

\newcommand\FN{F_\levelN} 
\newcommand\regr{\varrho} 
\DeclareMathOperator{\init}{in}
\begin{document}

\title{Generating Functions for Hecke Operators}

\author{Hala Hajj Shehadeh}
\address{Department of Mathematics, Courant Institute of Mathematical
  Sciences, 251 Mercer Street, New York, NY 10012, USA}
\email{hala@cims.nyu.edu}

\author{Samar Jaafar}
\address{Department of Mathematics, Rutgers University, Hill Center, 110
  Frelinghuysen Road, Piscataway, NJ 08854, USA}
\email{jaafar@math.rutgers.edu}

\author{Kamal Khuri-Makdisi}
\address{Corresponding author: Mathematics Department and Center for
Advanced Mathematical Sciences, American University of Beirut, Bliss
Street, Beirut, Lebanon}
\email{kmakdisi@aub.edu.lb}

\thanks{October 30, 2006}

\begin{abstract}
Fix a prime $N$, and consider the action of the Hecke operator $T_N$ on the
space $\mathcal{M}_\kappa(SL(2,\Z))$ of modular forms of full level and
varying weight $\kappa$.  The coefficients of the matrix of $T_N$ with
respect to the basis $\{E_4{}^i E_6{}^j \mid 4i + 6j = \kappa\}$ for
$\mathcal{M}_\kappa(SL(2,\Z))$ can be combined for varying $\kappa$ into a
generating function $F_N$.  We show that this generating function is a
rational function for all $N$, and present a systematic method for
computing $F_N$.  We carry out the computations for $N = 2, 3, 5$, and
indicate and discuss generalizations to spaces of modular forms of
arbitrary level.
\end{abstract}

\maketitle

\section{Introduction}
\label{section1}

In this article, we introduce a generating function for the action of a
Hecke operator on the spaces of elliptic modular forms of fixed level and
varying weight, and show that this generating function is a rational
function, whose coefficients belong to $\Q$ in many cases of interest.  We
start with the case of full level, i.e., of modular forms on the upper half
plane that are invariant (under the slash operator $f| \alpha$) with
respect to $SL(2,\Z)$.  In that setting, the graded ring $\ring$ of modular
forms on $SL(2,\Z)$ is generated by the Eisenstein series $E_4$ and $E_6$
of weights $4$ and~$6$, and we restrict ourselves for simplicity to the
Hecke operator $\TN$ with $\levelN$ prime.  Our generating 
function $\FN(a,b,A,B)$ is built up from the coefficients obtained when we
express each $\TN(E_4{}^i E_6{}^j)$ as a polynomial in $E_4$ and $E_6$.
(Actually, to simplify the calculations, we work instead with multiples of
$E_4$ and $E_6$ corresponding to the coefficients $a$ and $b$ in the
equation  $y^2 = x^3 + ax + b$ of an elliptic curve.)  Our main result is
then that $\FN$ is a rational function, i.e., a ratio of polynomials in the
variables $a,b,A,B$, with coefficients in $\Q$.  The rationality of $\FN$
generalizes to Hecke operators on modular forms of any level, and the
coefficients in the rational functions belong to $\Q$ in most cases of
interest; in general, however, the coefficients can belong to a cyclotomic
field.

The fact that $\FN$ and its generalizations are rational functions holds in
a rather general setting that extends beyond elliptic modular forms,
provided we have a finitely generated $\C$-algebra of modular forms of
fixed level and varying weight (i.e., varying type at infinity).  In some
sense, this generalizes the fact that the generating function $F_1$ for the
``identity'' Hecke operator $T_1$ is related to the Hilbert (or
Poincar\'e) series of the graded algebra of modular forms of fixed
level and varying weight.  Our strategy is to express the action of $\TN$
in general as a trace between two such graded algebras of modular forms,
extending the usual trace of modular forms from a smaller congruence
subgroup $\Gamma'$ to a larger subgroup $\Gamma$, defined by a sum over
cosets of $\Gamma' \backslash \Gamma$.

It would be interesting to try using the generating function $\FN$ and its
analogs to tackle some conjectures about $p$-adic slopes of Hecke
operators, such as for example Conjecture~1 of \cite{BuzzardCalegari} about
the $2$-adic slopes of $T_2$ acting on cusp forms for $SL(2,\Z)$.  Several
other authors have made computations of such slopes in the context of
overconvergent modular forms, and the appearance of rational generating
functions for Hecke operators has been noticed before~\cite{Smithline}, but
the rational generating functions there are computed \emph{ad hoc} from a
less conceptual perspective than our generalized Hilbert series of the
graded algebra of modular forms.  After this paper was completed, we also
learned of independent unpublished work by K. Buzzard (circa 1995) and
H. Hida (circa 1985) on some cases of generating functions for Hecke
operators. 

Although we have not pursued the matter in this article, we point out that
the rationality of $\FN$ and its generalizations also implies the
rationality of the generating function for the \emph{traces} of $\TN$.
This can be proved directly in a general setting from the trace formula, as
in~\cite{FrechetteOnoPapanikolas}.

Beyond giving an existence proof for the rationality of $\FN$, we describe
a framework that allows us to compute the rational function directly for
arbitrary prime $\levelN$ when we work with full level $SL(2,\Z)$.  The
calculations can be done without any reference to $q$-expansions, simply
via calculations of isogenies~\cite{Velu} from an elliptic curve 
$y^2 = x^3 + ax + b$ for transcendental $a, b$ to its quotient by a cyclic
subgroup of order $\levelN$.  Our approach allows us to express $\FN$ as
the trace of the inverse of a matrix with polynomial entries in $a,b,A,B$.
Something similar can in principle be done for a fixed Hecke 
operator on arbitrary level, provided we can compute convenient generators
and relations for the $\C$-algebra of all modular forms of that level.  In
any case, the computations involve fairly large numbers and matrices fairly
quickly, and we have contented ourselves with describing the generating
functions for $T_2$, $T_3$, and $T_5$ on the full modular group
$SL(2,\Z)$.

The senior author (KKM) would like to thank the junior authors (HHS and SJ) 
for joint work on understanding the formulas in~\cite{Velu} and in carrying
out the calculations in Sections \ref{section4} and~\ref{section5}.  We
also thank S. Abu Diab for some assistance with the calculations.  The
senior author also gratefully acknowledges support from the University
Research Board, at the American University of Beirut, and from the Lebanese
National Council for Scientific Research, for supporting this work through
the grants ``Equations for modular and Shimura curves.''

\section{Hecke operators and traces on algebras}
\label{section2}

Let $\Gamma \subset SL(2,\Z)$ be a congruence subgroup.  We write
$\modformsk(\Gamma)$ for the space of weight $\wtk$ modular forms which are
invariant with respect to $\Gamma$.  Viewing the forms as holomorphic
functions on the upper half plane $\Half$, we can multiply two modular
forms of weights $\wtk$ and $\wtk'$ to get a product of weight
$\wtk+\wtk'$.  Thus we can combine forms of all weights with respect to a
fixed group $\Gamma$ into a graded ring of modular forms
\begin{equation}
  \label{equation2.1}
\ring_\Gamma = \Directsum_{\wtk \geq 0} \modformsk(\Gamma).  
\end{equation}
We abbreviate $\ring = \ring_{SL(2,\Z)}$.

We use the usual notation $\Gamma_0(\levelN)$, $\Gamma_1(\levelN)$, and
$\Gamma(\levelN)$ for the congruence subgroups of $\Gamma(1) = SL(2,\Z)$.
We know that the ring $\ring$ of modular forms on $\Gamma(1)$ is generated
by the Eisenstein series of weights $4$ and~$6$, which we multiply below
by suitable constants.  We thus obtain modular forms $a(\tau)$ and
$b(\tau)$, for $\tau \in \Half$, such that 
\begin{equation}
  \label{equation2.3}
\begin{split}
\ring &= \ring_{\Gamma(1)} = \C[a,b], \\
a = a(\tau) &= -\frac{15}{\pi^4} \sum_{0 \neq \ell \in L_\tau} \ell^{-4}
   = -\frac{1}{3}
           \Bigl(1 + 240 \sum_{n = 1}^\infty \sigma_3(n) q^n\Bigr), \\
b = b(\tau) &= \frac{35}{\pi^6} \sum_{0 \neq \ell \in L_\tau} \ell^{-6}
  = \frac{2}{27}
           \Bigl(1 - 504 \sum_{n = 1}^\infty \sigma_5(n) q^n\Bigr).
\end{split}
\end{equation}
Here $L_\tau = \Z + \Z\tau$ is the lattice in $\C$ with basis
$\{1, \tau\}$, and $q = \exp(2\pi i \tau)$, as usual.  We have chosen the
above normalization to simplify the equation for the elliptic curve $E_\tau
= C/L_\tau$, as well as the analytic isomorphism between $\C/L_\tau$ and
$E_\tau$: 
\begin{equation}
  \label{equation2.4}
\begin{split}
  E_\tau :  y^2 &= x^3 + ax + b, \\
 P = z + L_\tau \in \C/L_\tau & \mapsto (x_P, y_P) = 
 \bigl(-\frac{\wp(z; L_\tau)}{\pi^2}, \frac{i \wp'(z; L_\tau)}{2\pi^3}\bigr)
                        \in E_\tau.
\end{split}
\end{equation}
In the above normalization, the invariant differential $\omega = dx/y$
on $E_\tau$ corresponds to $2\pi i\, dz$, where $z$ is the ``coordinate''
on $\C/L_\tau$.  We shall interchangeably consider modular forms as
\begin{itemize}
\item functions on $\Half$; 
\item functions of the triple $(E, \omega, \text{level structure})$ for an
  elliptic curve $E$ and of an invariant differential $\omega$ on $E$;
\item functions on pairs  $(L, \text{level structure})$ where $L \subset
  \C$ is a lattice.
\end{itemize}
We will omit mention of the level structure whenever this simplifies our
notation.  Thus $f(\tau) = f(E_\tau, \omega) = f(L_\tau)$ in our
three ways to view a modular form.

To simplify our treatment, we only consider the Hecke operator $\TN$ when
$\levelN$ is prime.  We use the usual normalization that is convenient for
$q$-expansions,
\begin{equation}
  \label{equation2.5}
f \in \modformsk(\Gamma(1)) \implies
  \TN f(\tau) = \frac{1}{\levelN}
               \Bigl(\levelN^\wtk f(\levelN\tau)
                     + \sum_{t = 0}^{\levelN-1} f((\tau + t)/\levelN)
               \Bigr).
\end{equation}
In terms of elliptic curves or lattices, this means that 
\begin{equation}
  \label{equation2.6}
\begin{split}
\TN f(E, \omega) &=
  \frac{1}{\levelN} \sum_{N\text{-isogenies } \pi: E\to E' \text{ with }
  \pi^*(\omega') = \omega} f(E', \omega'), \\
\TN f(L_\tau) &=
  \frac{1}{\levelN} \sum_{L \subset L' \text{ with } [L':L] = \levelN} f(L').
\end{split}
\end{equation}

We extend the Hecke operator $\TN$ additively, so that we have a
$\C$-linear map $\TN: \ring \to \ring$.  We also define a linear map
(depending on $\levelN$), sending $f \in \ring$ to $f' \in \ringoN$, by
additively extending the definition
\begin{equation}
  \label{equation2.7}
f \in \modformsk(\Gamma(1)) \implies f'(\tau) = \levelN^\wtk f(\levelN \tau).
\end{equation}
If $\Gamma'$ and $\Gamma$ are congruence subgroups with $\Gamma' \subset
\Gamma$, then we define a trace map
$\tr^{\Gamma'}_\Gamma: \modring(\Gamma') \to \modring(\Gamma)$ by
additively extending
\begin{equation}
  \label{equation2.8}
f \in \modformsk(\Gamma') \implies
   \tr^{\Gamma'}_\Gamma f = \sum_{\alpha \in \Gamma' \backslash \Gamma}
                               f |_\wtk \alpha. 
\end{equation}
Thus $\alpha$ ranges over a finite set of coset representatives giving a
disjoint union $\Gamma = \disjointunion_\alpha \Gamma' \alpha$.  As usual,
$(f |_\wtk \alpha) (\tau) = f(\alpha \tau) j(\alpha, \tau)^{-\wtk}$. 
We immediately obtain the following result: 

\begin{proposition}
  \label{proposition2.1}
\begin{enumerate}
\item
The map $f \in \ring \mapsto f' \in \ringoN$ is a homomorphism of algebras.
\item
The map $\TN: \ring \to \ring$ is given by 
\begin{equation}
  \label{equation2.9}
\TN f = \frac{1}{\levelN} \troN f',
\end{equation}
for the trace map $\troN: \ringoN \to \ring$.
\end{enumerate}
\end{proposition}

The ``matrix'' of $\TN$ with respect to the basis
$\{a^i b^j \mid i,j \geq 0\}$ of $\ring$ is described by coefficients
$c^\levelN_{ijkl}$, for $i,j,k,\ell \geq 0$, where
\begin{equation}
  \label{equation2.10}
\TN(a^i b^j) = \sum_{k,l \geq 0} c^\levelN_{ijkl} a^k b^l.
\end{equation}
(By looking at weights of forms, we see that $c^\levelN_{ijkl} = 0$ unless
$4i+6j = 4k+6l$, so the above sum is finite.)  We combine these
coefficients $c^\levelN_{ijkl}$ into a power series, thereby obtaining the
following generating function:
\begin{equation}
  \label{equation2.11}
\begin{split}
\FN(a,b,A,B) &=  \sum_{i,j,k,l \geq 0} c_{ijkl} A^i B^j a^k b^l \\ 
             &= \sum_{i,j \geq 0} A^i B^j \TN(a^i b^j)
                                \in \ring[[A,B]] = \C[a,b][[A,B]]. 
\end{split}
\end{equation}
Here $A,B$ are formal variables, and we can view $a,b$ as independent
transcendental variables as well.

\section{The case $\levelN = 2$}
\label{section3}

In this section, we have $\Gamma_0(2) = \Gamma_1(2)$, and the level
structure that this parametrizes on a given elliptic curve $E$ is a $2$-torsion
point $(e,0) \in E$, corresponding to $P = 1/2 + L_\tau \in \C/L_\tau$.
We can view $e = x_P$ as the weight~$2$ Eisenstein series 
\begin{equation}
  \label{equation3.1}
e(\tau) = -\pi^{-2} \wp(1/2; L_\tau)
  = - \frac{2}{3}
    \Bigl[1 + 24 \sum_{n\geq 1}
                         \Bigl(\sum_{\text{odd } d|n} d\Bigr)
                          q^n
    \Bigr]
\in \modforms_2(\Gamma_1(2)).
\end{equation}
(The above identity follows from the Fourier expansion of $\wp(z; L_\tau)$
in terms of $q$ and $\exp(2\pi i z)$.)
We of course have the identity $e^3 + ae + b = 0$ in $\ringo{2} =
\ringi{2}$.  It follows that the two forms $a$ and $e$ are
algebraically independent over $\C$, since otherwise $a$ and $b$ would be
algebraically dependent.

\begin{proposition}
  \label{proposition3.1}
The algebra $\ringo{2}$ is generated by $a$, $b$, and $e$ in weights $4$,
$6$, and $2$, respectively, subject only to the relation $e^3 + ae + b =
0$.  Thus we have
\begin{equation}
  \label{equation3.2}
\ringo{2} = \C[a,b,e]/(e^3 + ae + b) = \ring[e]/(e^3 + ae + b),
\end{equation}
and hence $\ringo{2}$ is a free $\ring$-module of rank~$3$, with basis
$\{1, e, e^2\}$.
\end{proposition}
\begin{proof}
We sketch a proof of this standard result. We have inclusions of graded
rings $\ring = \C[a,b] \subset \C[a,b,e]/(e^3 + ae + b) \subset \ringo{2}$.
The Hilbert series (with respect to the weight) of $\ringo{2}$ is
$1/(1-t^2)(1-t^4) = (1+t^2+t^4)/(1-t^4)(1-t^6)$ by standard formulas for
the dimension $\dim \modformsk(\Gamma_0(2))$.  On the other hand,
the subring $\C[a,b,e]/(e^3 + ae + b)$ already has the same Hilbert series
as all of $\ringo{2}$, so they must be equal.  Note that since
$b = -e^3 - ae$, we could have phrased our result as $\ringo{2} = \C[a,e]$;
however, we primarily wish to view $\ringo{2}$ as an $\ring$-module, for
the purpose of computing traces.  
\end{proof}

We can now represent elements of $\ringo{2}$ as $3\times 3$ matrices with
elements in $\ring$ by the regular representation of $\ringo{2}$ as an
$\ring$-algebra, with respect to, say, the basis $\{1, e, e^2\}$.  We thus
obtain:

\begin{corollary}
  \label{corollary3.2}
There exists a unique homomorphism of $\C$-algebras
$\regr : \ringo{2} \to M_{3\times 3}(\ring)$ such that
\begin{equation}
  \label{equation3.3}
\regr(a) =
   \begin{pmatrix} a & 0 & 0 \\ 0 & a & 0 \\ 0 & 0 & a \end{pmatrix}, 
\qquad
\regr(b) =
   \begin{pmatrix} b & 0 & 0 \\ 0 & b & 0 \\ 0 & 0 & b \end{pmatrix}, 
\qquad
\regr(e) =
    \begin{pmatrix} 0 & 0 & -b \\ 1 & 0 & -a \\ 0 & 1 & 0 \end{pmatrix}.
\end{equation}
The trace map from $\ringo{2}$ to $\ring$ can be computed as the
matrix trace:
\begin{equation}
  \label{equation3.4}
f \in \ringo{2} \implies \tr^{\Gamma_0(2)}_{\Gamma(1)} f = \tr \regr(f).
\end{equation}
\end{corollary}
\begin{proof}
We introduce the fields of fractions $\mathcal{K}$,
$\mathcal{K}_{\Gamma_0(2)}$, and $\mathcal{K}_{\Gamma(2)}$ 
of the integral domains $\ring$, $\ringo{2}$, and $\ring_{\Gamma(2)}$;
since we have taken fields of fractions without respecting the graded
structure, these fields are less natural than function fields of modular
curves or than graded rings of \emph{meromorphic} modular forms (obtained
by inverting only the nonzero homogeneous elements).  The action $f \mapsto
f|\alpha$ gives the integral domain $\ring_{\Gamma(2)}$ an action 
of the finite group $G = \Gamma(1)/\Gamma(2) \isomorphic SL(2, \Z/2\Z)$,
such that $\ring$ is the subring invariant under $G$, and $\ringo{2}$ is
the subring invariant under $H = \Gamma_0(2)/\Gamma(2) \subset G$.  This
situation is mirrored in the fields of fractions, the main point being that
elements of $\mathcal{K}_{\Gamma(2)}$ can be written in the form $f/s$ with
a $G$-invariant denominator $s \in \ring$.  (Given an arbitrary
denominator, multiply above and below by the ``conjugates'' of the
denominator.)  This last observation, combined with
Proposition~\ref{proposition3.1}, also shows that $\{1,e,e^2\}$ is a basis
for $\mathcal{K}_{\Gamma_0(2)}$ as a vector space over $\mathcal{K}$.
Elementary Galois theory applied to $\mathcal{K}$,
$\mathcal{K}_{\Gamma_0(2)}$, and $\mathcal{K}_{\Gamma(2)}$ now yields our
statement that the trace of the regular representation is equal to the sum
over conjugates given by the action of representatives $\alpha$ for
$H\backslash G \isomorphic \Gamma_0(2)\backslash \Gamma(1)$.  
\end{proof}

\begin{proposition}
  \label{proposition3.3}
The ring homomorphism $f \in \ring \mapsto f' \in \ringo{2}$ of
\eqref{equation2.7} satisfies
\begin{equation}
  \label{equation3.5}
a' = -4a - 15e^2, \qquad\qquad b' = 22b + 14ae.  
\end{equation}
\end{proposition}
\begin{proof}
We know that $a' \in \modforms_4(\Gamma_0(2)) = \C a + \C e^2$ and
$b' \in \modforms_6(\Gamma_0(2)) = \C b + \C ae$.  Thus we must find
constants $\mu_1, \dots, \mu_4$ such that $a' = \mu_1 a + \mu_2 e^2$ and
$b' = \mu_3 b + \mu_4 ae$.  We can do this in one of three ways: the first
way is to compare $q$-expansions, the second way is to use the known traces
of the Hecke operator $T_2$ on the spaces $\modformsk(\Gamma(1))$ for $\wtk
\leq 12$, and the third way is to use  V\'elu's formulas for isogenies, as
described in Section~\ref{section4} below.
\end{proof}

\begin{theorem}
  \label{theorem3.4}
The generating function $F_2(a,b,A,B)$ of \eqref{equation2.11} is
\begin{equation}
  \label{equation3.6}
F_2(a,b,A,B)
  = \frac{1}{2}
    \tr \Bigl[
       (I - A \regr(a')) (I - B \regr(b'))
        \Bigr]^{-1} 
  = \frac{1}{2} \tr M^{-1},
\end{equation}
where $I$ is the $3\times3$ identity matrix, and $M$ is the product
\begin{equation}
\label{equation3.6.5}
  M =  
     \begin{pmatrix}
            1+4Aa  & -15Ab  & 0      \\
            0      & 1-11Aa & -15Ab  \\
            15A    & 0      & 1-11Aa
     \end{pmatrix}
     \begin{pmatrix}
            1-22Bb & 0      & 14Bab  \\
            -14Ba  & 1-22Bb & 14Ba^2 \\
            0      & -14Ba  & 1-22Bb
     \end{pmatrix}.
\end{equation}
\end{theorem}
\begin{proof}
We extend the operations $f\mapsto f'$ and $\tr^{\Gamma_0(2)}_{\Gamma(1)}$
coefficientwise, so that they map from $\ring[[A,B]]$ to $\ringo{2}[[A,B]]$
and vice-versa.  We similarly extend the regular representation $\regr$ so
that it sends $\ringo{2}[[A,B]]$ to $M_{3\times 3}(\ring[[A,B]])$.  All
relations are still valid after this extension of scalars; in particular
$\tr^{\Gamma_0(2)}_{\Gamma(1)}$ can again be computed as a matrix trace.
Now by the results of Section~\ref{section2}, and by linearity in each
coefficient of a monomial $A^i B^j$, we have
\begin{equation}
  \label{equation3.7}
F_2(a,b,A,B)
  = \frac{1}{2}
    \tr^{\Gamma_0(2)}_{\Gamma(1)}
     \Bigl[
      \sum_{i,j \geq 0}
         A^i B^j (a')^i (b')^j
    \Bigr]
 =  \frac{1}{2}
    \tr^{\Gamma_0(2)}_{\Gamma(1)}
     \left[
          \frac{1}{(1-Aa')(1-Bb')}
    \right].
\end{equation}
We know the values for $a',b'$ from~\eqref{equation3.5}.  We then compute
the trace as in~\eqref{equation3.4}; note that the regular representation
$\regr$ respects inverses in the power series ring $\ring[[A,B]]$.
\end{proof}

Note that it is clear from the above that $F_2$ is a rational function in
$\Q(a,b,A,B)$.  Moreover, the denominator of $F_2$ is the determinant
$\det M$, which appears in the denominator of $M^{-1}$.

\section{General framework for $\TN$ in terms of V\'elu's formulas}
\label{section4}

In this section, $\levelN$ is an odd prime.  One can prove a version of the
results below even if $\levelN=2$, or indeed if $\levelN$ is composite, but
the statements and proofs become more complicated.

By the same reasoning as in Theorem~\ref{theorem3.4}, we deduce that
\begin{equation}
  \label{equation4.1}
\FN(a,b,A,B)
   = \frac{1}{\levelN} \troN \left[ \frac{1}{(1-Aa')(1-Bb')} \right].
\end{equation}
At this point, $a',b' \in \ringoN$, and we wish to find their values and be
able to represent them as matrices via a suitable generalization of the
regular representation $\regr$ from Section~\ref{section3}.
It is enough to work on the level of the fields of fractions $\mathcal{K},
\mathcal{K}_{\GoN}$ introduced in the proof of
Corollary~\ref{corollary3.2}, and to calculate with the regular 
representation with respect to a $\mathcal{K}$-basis of
$\mathcal{K}_{\GoN}$.  We introduce a further simplification by passing to
the smaller congruence subgroup 
\begin{equation}
  \label{equation4.2}
\GpmN = \{\pm I \} \cdot \GiN
    \subset \GoN,
\end{equation}
for which, as we shall see, the corresponding field of fractions
$\mathcal{K}_{\GpmN}$ is easier to describe and to study via a regular
representation than the original $\mathcal{K}_{\GoN}$.  Thus we can easily
evaluate the desired trace with respect to this smaller group.  The
resulting trace is off by a factor $[\GpmN:\GoN] = (\levelN - 1)/2$, so we
retrieve the original trace by dividing:
\begin{equation}
  \label{equation4.3}
 \troN [(1-Aa')(1-Bb')]^{-1}
  = \frac{2}{\levelN-1}
    \tr^{\GpmN}_{\Gamma(1)} [(1-Aa')(1-Bb')]^{-1}.
\end{equation}

In order to describe $\mathcal{K}_{\GpmN}$, we begin with the level
structure parametrized by $\GiN$, namely the $\levelN$-torsion point on our
varying elliptic curve, corresponding to $P = 1/\levelN + L_\tau \in
\C/L_\tau$.  Since $\GpmN$ introduces an ambiguity between $P$ and $-P$,
the $x$-coordinate $x_P$ of $P$ is invariant under $\GpmN$.  Hence we
obtain 
\begin{equation}
  \label{equation4.4}
x_P(\tau) = - \pi^{-2} \wp(1/\levelN;L_\tau) \in \modforms_2(\GpmN),
\end{equation}
and, more generally,
\begin{equation}
  \label{equation4.5}
x_{\ell P}(\tau) = - \pi^{-2} \wp(\ell/\levelN;L_\tau) 
                                      \in \modforms_2(\GpmN), 
\quad \ell \in S = \{1, \dots, (N-1)/2\}.
\end{equation}
(The above modular forms are, incidentally, all Eisenstein series.)  We do
not need other values of $\ell$, since $x_{-\ell P} = x_{\ell P}$.  We
begin with the following proposition:

\begin{proposition}
  \label{proposition4.1}
Let $\psi_\levelN(x;a,b) = \levelN x^{(\levelN^2-1)/2} + \dots \in \Z[a,b,x]$
be the $N$-division polynomial (see, e.g., Exercise III.3.7
of~\cite{SilvermanI}).
\begin{enumerate}
\item
We have $\mathcal{K}_{\GpmN} = \mathcal{K}[x_P]/(\psi_\levelN(x_P;a,b))$.
\item
The powers $\{1, x_P, x_P{}^2, \dots, x_P{}^{(\levelN^2-3)/2}\}$
are a basis for $\mathcal{K}_{\GpmN}$ over $\mathcal{K}$.
\item
The other modular forms $x_{\ell P}$ for $\ell \in S$ also belong to
$\mathcal{K}_{\GpmN}$; their expressions in terms of $x_P$, $a$, and
$b$ are straightforward to compute, and involve only coefficients from
$\Q$.
\item
Every symmetric polynomial in the $x_{\ell P}$ belongs to $\ringoN$.
\end{enumerate}
\end{proposition}
\begin{proof}
By Galois theory for the extension
$\mathcal{K}_{\Gamma(\levelN)}/\mathcal{K}$, we know that $x_P$ generates
$\mathcal{K}_{\GpmN}$ over $\mathcal{K}$, because $x_P$ is left invariant
precisely by the subgroup $\GpmN$ of $\Gamma(1)$.  Moreover, $x_P$ is a root
of the division polynomial $\psi_\levelN(x;a,b)$, which is irreducible over
$\mathcal{K}$ because $\Gamma(1)$ acts transitively on the (nonzero)
torsion points in $E[\levelN]$, and hence on the roots of $\psi_\levelN$.
(Note that if $\levelN$ were not prime, we would need to work with the
``primitive'' $\levelN$-division polynomial instead of $\psi_\levelN$.)
This proves the first two statements of the proposition.  The third
statement holds from the multiplication formula $x_{\ell P} =
\phi_\ell(x_P;a,b)/[\psi_\ell(x_P;a,b)]^2$, where $\psi_\ell$ is
the $\ell$-division polynomial, and $\phi_\ell \in \Z[a,b,x]$ is another
polynomial that is straightforward to compute.  The fourth statement is
immediate, since the action of an element of $\GoN$ transforms $P$ into a
multiple $\ell' P$ with $(\ell', \levelN) = 1$.
\end{proof}

In light of the above discussion, it is straightforward to calculate $\FN$
provided we can express $a'$ and $b'$ as elements of
$\mathcal{K}[x_P]/(\psi_\levelN(x_P;a,b))$.  This is not as easily done as
in the case $\levelN=2$, where we simply compared $q$-expansions, since the
expressions that we seek for $a'$ and $b'$ are rational functions and not
necessarily polynomials in $x_P$, $a$, and $b$.  Thus, if we wished to use
$q$-expansions in order to find expressions for $a', b'$, we would need to
bound the denominators of those expressions.  We instead compute $a', b'$
using the interpretation of modular forms as functions of elliptic curves
$E$ with a choice of global differential $\omega$.

\begin{proposition}
  \label{proposition4.2}
We have
\begin{equation}
  \label{equation4.6}
\begin{split}
  a' &= 
    a - 30 \sum_{\ell \in S} x_{\ell P}{}^2
      - 5(\levelN-1) a,
\\
  b' &=
    b - 70 \sum_{\ell \in S} x_{\ell P}{}^3 
      - 42a \sum_{\ell \in S} x_{\ell P}
      - 14(\levelN-1) b.
\end{split} 
\end{equation}
\end{proposition}
\begin{proof}
According to V\'elu's formulas for isogenies~\cite{Velu}, the coefficients
$a', b'$ above are the coefficients in a Weierstrass equation
$E': {y'}^2 = {x'}^3 + a'x' + b'$ for the quotient curve
$E' = E/\langle x_P \rangle$ of the elliptic curve $E: y^2 = x^3 + ax + b$
by the subgroup generated by $P \in E[\levelN]$.  (Recall that $\levelN$
is odd; V\'elu's formulas are slightly different when $2$-torsion is
involved.)  The isogeny given by the projection $\pi: E \to E'$ is such
that the pullback of the global differential $\omega' = dx'/y'$ 
is $\pi^* \omega' = \omega = dx/y$.  Thus $a'$ (respectively, $b'$)
corresponds to one term in the sum over isogenies defining $\TN a$
(respectively, $\TN b$) in~\eqref{equation2.6}.  In fact, $a'$
specifically corresponds to the term $\levelN^\wtk a(\levelN \tau)$
in~\eqref{equation2.5}, and similarly for $b'$.  The reason is that if we
view $E$ as $\C/L_\tau$, $\omega$ as $2\pi i \,dz$, and $P$ as the image of
$1/\levelN$, then $E'$ is $\C/(\Z\frac{1}{\levelN} + \Z\tau)$ with
$\omega'=2 \pi i \, dz$ as well;
hence $a'(\tau) = a(\Z\frac{1}{\levelN} + \Z\tau) 
= \levelN^4 a(L_{\levelN\tau})$, in the interpretation of $a = a(L)$ as a
function of lattices.
\end{proof}

Putting together the results of this section, we obtain the main result of
this article:

\begin{theorem}
  \label{theorem4.3}
The generating function $\FN$ defined in~\eqref{equation2.11} is a rational
function in $\Q(a,b,A,B)$, and $\FN$ can be explicitly computed for any
specific value of $\levelN$.
\end{theorem}
\begin{proof}
The case $\levelN = 2$ is Theorem~\ref{theorem3.4}.  For $\levelN \geq 3$,
combine equations \eqref{equation4.1} and~\eqref{equation4.3}, as well as
Propositions \ref{proposition4.1} and~\ref{proposition4.2}.  The trace
$\tr^{\GpmN}_{\Gamma(1)}$ in~\eqref{equation4.3} can be computed via the
regular representation of $\mathcal{K}_{\GpmN}$ as a $\mathcal{K}$-algebra,
with respect to the basis consisting of the powers of $x_P$.  All the
constants that we encounter belong to $\Q$, most significantly by the third
statement in Proposition~\ref{proposition4.1}.
\end{proof}

We remark that it is unfortunate that we need to use square matrices of size
$(\levelN^2-1)/2 = [\Gamma(1):\GpmN]$ in our regular representation of
$\mathcal{K}_{\GpmN}$, since the true dimension that matters is $\levelN+1
= [\Gamma(1):\GoN]$.  It would be agreeable to have a direct way to
describe a cyclic $\levelN$-subgroup of $E$ and the isogeny obtained by
quotienting $E$ by that subgroup.  This would be simpler than our 
approach of choosing the $N$-torsion point $P$ and considering the subgroup
$\langle P \rangle$ generated by $P$.

\section{Calculations for $T_3$ and $T_5$}
\label{section5}

In this section, we give explicit matrices $\regr(a')$ and $\regr(b')$, in
the two cases $\levelN = 3$ and $\levelN = 5$, for a suitable regular
representation $\regr$ of $\mathcal{K}_{\GoN}$.  (Note that $\Gamma_0(3) =
\Gamma_{\pm}(3)$, but $\Gamma_0(5) \neq \Gamma_{\pm}(5)$, so we must modify
the approach in Section~\ref{section4} when $\levelN=5$.)  This is enough
data to describe $F_3$ and $F_5$ completely, by evaluating the trace
in~\eqref{equation4.1} using $\regr$.

As noted in the above paragraph, the case $\levelN=3$ is exactly covered by
our previous methods.  In this case $\psi_3 = 3x^4 + 6ax^2 + 12bx - a^2$,
so $\regr(a) = aI$ and $\regr(b) = bI$, where $I$ is the $4\times 4$
identity matrix, and
\begin{equation}
  \label{equation5.1}
  \regr(x_P)
  = \begin{pmatrix}
       0 & 0 & 0 & a^2/3 \\
       1 & 0 & 0 & -4b   \\
       0 & 1 & 0 & -2a   \\
       0 & 0 & 1 & 0
    \end{pmatrix}.
\end{equation}
We fortunately have $S = \{1\}$, so $a' = -9a-30x_P{}^2$ and $b' = -27b -
70x_P{}^3 - 42ax_P$.  Hence we obtain the following result.
\begin{theorem}
  \label{theorem5.1}
For $\levelN = 3$, we have that $\regr(a')$ and $\regr(b')$ are respectively
\begin{equation}
  \label{equation5.2}
\begin{pmatrix}
  -9a & 0   & -10a^2 & 0      \\
    0 & -9a & 120b   & -10a^2 \\
  -30 & 0   & 51a    & 120b   \\
    0 & -30 & 0      & 51a    \\
\end{pmatrix}
\text{ and }
\begin{pmatrix}
 -27b & -70a^2/3 & 0        & 98a^3/3   \\
 -42a & 253b     & -70a^2/3 & -392ab    \\
    0 & 98a      & 253b     & -658a^2/3 \\
  -70 & 0        & 98a      & 253b      \\
\end{pmatrix}.
\end{equation}
\end{theorem}

We now turn to the case $\levelN=5$.  In this setting, we first worked
as described in Section~\ref{section4}, using the regular representation of
$\mathcal{K}_{\Gamma_{\pm}(5)}$ with respect to the basis
$\{1, x_P, x_P{}^2, \dots, x_P{}^{11}\}$.  We were however dissatisfied
with the appearance of the results --- our first calculations gave us
$12 \times 12$ matrices $\regr(a'), \regr(b')$ whose entries all had a
denominator of $(4a^3 + 27b^2)^2$, i.e., the square of the discriminant.
We preferred instead to work with elements invariant under $\Gamma_0(5)$,
and so looked for ``nice'' symmetric polynomials in $x_P, x_{2P}$
(corresponding to $S=\{1,2\}$) that gave a basis for the $6$-dimensional
field extension $\mathcal{K}_{\Gamma_0(5)}/\mathcal{K}$.  After some trial
and error, we settled on the following ordered basis, which gave us
matrices with polynomial entries:
\begin{equation}
  \label{equation5.3}
\{f_1, \dots, f_6\} =
\{1, x_P + x_{2P}, x_P x_{2P}, x_P{}^2 + x_{2P}{}^2,
          x_P{}^2 x_{2P}{}^2, x_P{}^3 + x_{2P}{}^3 \}.
\end{equation}

\begin{remark}
\label{remark5.2}
To help the reader check the calculations, we mention that in the regular
representation $\regr$ with respect to the above basis $\{f_1, \dots,
f_6\}$, we have 
\begin{equation}
  \label{equation5.4}
\regr(x_P+x_{2P}) = 
\begin{pmatrix}
0 & 0 & -4b/3 & -4b/3 & 16ab/3  & 3a^2 \\
1 & 0 & -2a/3 & -2a/3 & 11a^2/3 & -8b  \\
0 & 2 & 0     & 0     & -12b    & -2a  \\
0 & 1 & 0     & 0     & -4b     & 2a   \\
0 & 0 & 0     & 0     & 0       & 15   \\
0 & 0 & 1/3   & 4/3   & -4a/3   & 0  
\end{pmatrix},
\end{equation}
\begin{equation}
  \label{equation5.5}
\regr(x_P x_{2P}) = 
\begin{pmatrix}
0 & -4b/3 & 0 & a^2 & (16b^2-4a^3)/5 & 40ab/3  \\
0 & -2a/3 & 0 & -4b & 24ab/5         & 29a^2/3 \\
1 & 0     & 0 & -2a & 9a^2/5         & -32b    \\
0 & 0     & 0 & 0   & a^2/5          & -12b    \\
0 & 0     & 1 & 3   & -18a/5         & 0       \\
0 & 1/3   & 0 & 0   & -8b/5          & -10a/3
\end{pmatrix}.
\end{equation}
This can be checked by working in the full field
$\mathcal{K}_{\Gamma_{\pm}(5)}$ in terms of the basis of powers of $x_P$.
Alternatively, the reader may wish to verify the above matrices by working
directly from the algebraic relations satisfied by $x_P$ and $x_{2P}$ over
$\mathcal{K}$.  Instead of using the cumbersome fact that $x_P$ and
$x_{2P}$ are roots of the high degree division polynomial $\psi_5$, it is
easier to note that $x_{4P} = x_P$ and to use the duplication formula for
points on $E$ to deduce the relations
$[\psi_2(x_P)]^2 x_{2P} = \phi_2(x_P)$ and
$[\psi_2(x_{2P})]^2 x_P = \phi_2(x_{2P})$.
(We also need the fact that $x_P \neq x_{2P}$ since $P$ does not have order
$3$; the ideal of relations between $x_P$ and $x_{2P}$ can be obtained
by starting with the ideal generated by the two formulas above and by
saturating that ideal with respect to $x_P - x_{2P}$.)
\end{remark}

Using \eqref{equation5.4} and~\eqref{equation5.5}, we now easily find
$\regr(x_P{}^2 + x_{2P}{}^2)$ and $\regr(x_P{}^3 + x_{2P}{}^3)$, which
allow us to apply V\'elu's formulas to obtain the following result:
\begin{theorem}
\label{theorem5.3}
In the case $\levelN=5$, let $\regr$ be the regular representation of
$\mathcal{K}_{\Gamma_0(5)}$ over $\mathcal{K}$ with respect to the basis 
$\{f_1, \dots, f_6\}$ above.  Then we have
\begin{equation}
  \label{equation5.6}
\regr(a') = 
\begin{pmatrix}
-19a & 40b & -30a^2 & -60a^2 & 72a^3 - 448b^2 & -1600ab  \\
   0 & a   & 120b   & 120b   & -512ab         & -1160a^2 \\
   0 & 0   & 41a    & 0      & -192a^2        & 3960b    \\
 -30 & 0   & 0      & -79a   & -18a^2         & 1320b    \\
   0 & 0   & -90    & -420   & 365 a          & 0        \\
   0 & -40 & 0      & 0      & 184 b          & 321a
\end{pmatrix},
\end{equation}
\begin{equation}
  \label{equation5.7}
\regr(b') =
\begin{pmatrix}
-55b & -210a^2 & -2632ab/3  & -11032ab/3 & 13888a^2b/3  & 1554a^3-12320b^2 \\
-42a & 505b    & -1946a^2/3 & -8036a^2/3 & 7658a^3/3-3360b^2 & -13104ab    \\
0    & 56a     & 2185b      & 9240b      & -9576ab           & -5096a^2    \\
0    & -182a   & 840b       & 3025b      & -2632ab           & -994a^2     \\
0    & -1050   & 0          & 0          & 4705b             & 7770a       \\
-70  & 0       & 658a/3     & 2212a/3    & -2842a^2/3        & 5265b
\end{pmatrix}.
\end{equation}
\end{theorem}

\section{Generalizations}
\label{section6}

\subsection{Other levels than $\Gamma(1)$}
\label{subsection6.1}

Our first generalization is to study generating functions
for Hecke operators on $\ring_\Gamma$, for an arbitrary congruence subgroup
$\Gamma$.  Our approach can deal with any Hecke operator given by a double
coset $\Gamma \alpha \Gamma$ with $\alpha \in GL(2, \Q)$, 
$\det \alpha > 0$.  In this situation, we shall show in this subsection
that the analog of $\FN$ is still a rational function, with coefficients
in a number field; in many cases of interest, the coefficients actually lie
in $\Q$.  We can compute the analog of $\FN$ in any specific case, but we
do not have a satisfactory systematic method for computing the generating
function by methods analogous to those in Section~\ref{section4}.

The main issue in generalizing our previous argument to arbitrary $\Gamma$
is that as soon as the modular curve associated to $\Gamma$ has positive
genus, the ring $\ring_\Gamma$ is no longer a polynomial algebra in two
variables.  Hence the $\C$-basis $\{ a^i b^j \mid i, j \geq 0\}$ of
$\ring = \ring_{\Gamma(1)}$ that we used to define the coefficients
$c^\levelN_{ijkl}$ of~\eqref{equation2.10} must be replaced by something
more complicated for $\ring_\Gamma$.  We must do this in a way that still
yields an analog of the identity
$\sum_{i,j \geq 0} A^i B^j a^i b^j = [(1-Aa)(1-Bb)]^{-1}$
in $\ring[[A,B]]$ that plays such a crucial role in Theorems
\ref{theorem3.4} and~\ref{theorem4.3}.  We thus replace $a, b \in \ring$ by
generators $a_1, \dots, a_r \in \ring_\Gamma$, where $a_i \in
\modforms_{\wtk_i}(\Gamma)$.  (The ring $\ring_\Gamma$ is finitely
generated as a $\C$-algebra because, e.g., it us an integral extension
of $\ring$.)  Writing $I$ for the ideal of relations among the $\{a_i\}$,
we see that we need to find an appropriate $\C$-basis for $\ring_\Gamma =
\C[a_1, \dots, a_r]/I$.

\begin{proposition}
  \label{proposition6.1}
Given any finitely generated $\C$-algebra $\C[a_1, \dots, a_r]/I$, let
$\init(I)$ be the initial ideal of $I$ with respect to any fixed term order
on the monomials in the $\{a_i\}$.  Then the set
\begin{equation}
  \label{equation6.1}
\mathcal{B} = \{\text{monomials } m \mid m \notin \init(I)\}
\end{equation}
is a $\C$-basis for $\C[a_1, \dots, a_r]/I$, and the formal power series
\begin{equation}
  \label{equation6.2}
G(a_1, \dots, a_r) = \sum_{m \in \mathcal{B}} m \in \C[[a_1, \dots, a_r]]
\end{equation}
is actually a rational function, of the form
\begin{equation}
  \label{equation6.3}
G(a_1, \dots, a_r) = \frac{N(a_1, \dots, a_r)}{(1-a_1)(1-a_2)\dots(1-a_r)}
\end{equation}
where $N(a_1, \dots, a_r) \in \Z[a_1, \dots, a_r]$ is a polynomial with
integer coefficients.
\end{proposition}
\begin{proof}
The first assertion, that $\mathcal{B}$ is a basis, is a standard result in
the theory of Gr\"obner bases.  The second assertion, that $G$ is a
rational function with known denominator, follows from a direct
modification of the usual argument by induction on $r$ to show the
rationality of Hilbert 
series of graded modules.  We
apply this specifically to the module $M = \C[a_1, \dots, a_r]/\init(I)$,
which carries an action of the algebraic torus $T = (\C^*)^r$ such that an
element $t = (\lambda_1, \dots, \lambda_r) \in T$ sends $a_i$ to
$t a_i = \lambda_i a_i$, for $1 \leq i \leq r$.  Thus the series $G$ is the
same as the $T$-equivariant Hilbert series of $M$ discussed, e.g., in
Section~6.6 of~\cite{ChrissGinzburg}.
\end{proof}

\begin{remark}
  \label{remark6.2}
The function $G$ depends significantly on the choice of term order, to say
nothing of the choice of generators $a_1, \dots, a_r$.  Take for example
$\ringo{2} = \C[a,b,e]/(e^3 + ae + b)$.  Depending on whether the 
initial term in $e^3 + ae + b$ is $e^3$, $ae$, or $b$, we obtain
$G = (1+e+e^2)/[(1-a)(1-b)]$,
$G = (1-ae)/[(1-a)(1-b)(1-e)]$, or
$G = 1/[(1-a)(1-e)]$,  
respectively.  These examples incidentally show us that cancellation can
occur between the numerator and denominator of $G$.
\end{remark}

In light of the above proposition, we now see that the analog of our
earlier sum $\sum_{i,j} A^i B^j a^i b^j$ is
\begin{equation}
  \label{equation6.4}
G(A_1 a_1, \dots, A_r a_r) =
 \sum_{\text{ monomials } m = a_1^{i_1} \dots a_r^{i_r} \in \mathcal{B}} 
     A_1^{i_1} \dots A_r^{i_r} a_1^{i_1} \dots a_r^{i_r},
\end{equation}
while our analog of $\FN$, corresponding to the Hecke operator
$\Gamma\alpha\Gamma$, is 
\begin{equation}
  \label{equation6.5}
F_{\Gamma\alpha\Gamma}(A_1, \dots, A_r, a_1, \dots, a_r) =
 \sum_{m = a_1^{i_1} \dots a_r^{i_r} \in \mathcal{B}} 
     A_1^{i_1} \dots A_r^{i_r} 
     \Bigl[a_1^{i_1} \dots a_r^{i_r}\Bigm| \Gamma\alpha\Gamma \Bigr].
\end{equation}
We view the formal sums above in the ring
$\ring_\Gamma[[A_1, \dots, A_r]]$.
Our argument for the rationality of $F_{\Gamma\alpha\Gamma}$ now
proceeds essentially identically to our previous discussion, and we obtain
(up to normalization constants) an identity of the form
\begin{equation}
  \label{equation6.6}
F_{\Gamma\alpha\Gamma} = \tr^{\Gamma'}_\Gamma G(A_1 a'_1, \dots A_r a'_r), 
\quad \text{where }
\Gamma' = \Gamma \intersect \alpha^{-1} \Gamma \alpha
\text{ and } f' = f \mid \alpha.
\end{equation}
More precisely, depending on how we want to normalize the action of
$\Gamma\alpha\Gamma$, we can modify the definition of $f'$ by choosing a
constant $C$ and defining $f \in \modformsk(\Gamma) \implies f' =
C^\wtk f|_\wtk \alpha$; this ensures that the map $f \mapsto f'$ is still a
ring homomorphism from $\ring_\Gamma$ to $\ring_{\Gamma'}$.  We can also
include another constant factor in front of the trace
in~\eqref{equation6.6}.

The above suffices to show that $F_{\Gamma\alpha\Gamma}$ is a rational
function of the $\{a_i\}$ and the $\{A_i\}$, viewing these as independent
indeterminates.  The coefficients of this rational function can be taken to
lie in a field containing essentially the coefficients of the
$q$-expansions of all modular forms in $\ring_\Gamma$ and $\ring_{\Gamma'}$
that we encounter.

\begin{theorem}
  \label{theorem6.3}
Let $L$ be a subfield of $\C$ such that for all weights~$\wtk$, the spaces
$\modformsk(\Gamma)$ and $\modformsk(\Gamma')$ have a basis of forms whose
$q$-expansions have $L$-rational coefficients, and such that if $f \in
\modformsk(\Gamma)$ has $L$-rational coefficients, then so does $f'$.
Assume that the generators $\{a_i\}$ of $\ring_\Gamma$ are moreover chosen
to all have $L$-rational coefficients.  Then the generating function
$F_{\Gamma\alpha\Gamma}$ is a rational function in the indeterminates
$a_1, \dots, a_r, A_1, \dots, A_r$, with coefficients belonging to $L$.
\end{theorem}
\begin{proof}
The assumptions on the coefficients in the $q$-expansions allow us to
choose a $\mathcal{K}_\Gamma$-basis for $\mathcal{K}_{\Gamma'}$, with
respect to which the regular representation $\regr$ takes a form
$f \in \modformsk(\Gamma')$ with $L$-rational coefficients to a matrix
$\regr(f)$ whose entries are rational functions of the $\{a_i\}$ with
coefficients in $L$.  (Note that the ideal $I$ of relations between the
$\{a_i\}$ is also defined over $L$.)  Since the forms $\{a'_i\}$ are also
$L$-rational by assumption, we obtain our desired result.
\end{proof}

\begin{corollary}
  \label{corollary3.4}
We can always take $L$ above to be a cyclotomic field (provided the
normalizing constants like $C$ above also belong to $L$).  In the typical
case where $\Gamma$ is one of $\Gamma_0(\levelN)$, $\Gamma_1(\levelN)$, or
$\Gamma(\levelN)$, and $\alpha$ is a diagonal matrix, then we can take $L =
\Q$.  This typical case includes the ``standard'' Hecke operators $T_\ell$
on any one of the groups above, even when $(\ell, \levelN) > 1$.
\end{corollary}
\begin{proof}
It is well known that the space $\modformsk(\Gamma(\levelN))$ (similarly
for $\Gamma_0(\levelN)$ and $\Gamma_1(\levelN)$)
has a basis of forms whose $q$-expansions have $\Q$-rational coefficients,
but that replacing $f$ by $f|_{\wtk} \alpha$ where $\alpha$ is an integral
matrix with determinant $\ell$ can introduce roots of unity up to
$\levelN\ell$ (see for example Chapters 3 and~6 of~\cite{Shimura}).  Thus
the first assertion of our corollary is clear.  The second assertion
follows because the action of a diagonal matrix $\alpha$ replaces $\tau$ by
a multiple $a\tau/d$ for some $a,d \in \Z$, thereby preserving rationality
of $q$-expansions.
\end{proof}

\subsection{Hecke operators restricted to cusp forms}
\label{subsection6.2}

Our second generalization is that our method extends to construct rational
generating functions for Hecke operators acting only on the cuspidal part
$\mathcal{I}_\Gamma = \Directsum_{\wtk \geq 0} \mathcal{S}_{\wtk}(\Gamma)$
of the full ring of modular forms $\ring_\Gamma$.  This requires very
little work in light of our observations in
Subsection~\ref{subsection6.1}.  We merely need to point out that
$\mathcal{I}_\Gamma$ is an ideal in $\ring_\Gamma$.  Hence we can use a
Gr\"obner basis argument as in Proposition~\ref{proposition6.1} to produce
a rational function, analogous to $G$, which is the formal sum of a basis
for all cusp forms.  It is probably best to choose the generators $\{a_i\}$
of the $\C$-algebra $\ring_\Gamma$ to consist of cusp forms and Eisenstein
series, in such a way that $\mathcal{I}_\Gamma$ is generated by
$\{a_i \mid a_i \text{ is a cusp form}\}$.

\subsection{Automorphic forms on other groups}
\label{subsection6.3}

We conclude with the observation that the results in this article
generalize to other settings where one has graded rings of automorphic
forms.  This includes groups with Hermitian symmetric spaces, for which we
can interpret the automorphic forms as holomorphic functions on domains in
$\C^n$, such as the case of Hilbert modular forms (over a totally real
number field) and Siegel modular forms.  However, in the case of Hilbert
modular forms, we would probably be restricted to parallel weights, in
order to obtain a graded ring of automorphic forms that is a finitely
generated $\C$-algebra.  We do not need modular varieties with cusps to
carry out our program; modular forms on indefinite quaternion algebras over
$\Q$ come to mind, corresponding to automorphic forms on Shimura curves,
but it is more complicated to compute relations between the analogs of the
forms $a_i$ and $a'_i$ in that setting.

An alternative source of graded rings of automorphic forms is groups $G$
for which $G(\R)$ is compact, as discussed in~\cite{KKMRingsAutForms}.  In
that setting, one can view modular forms as holomorphic sections of line
bundles on several disjoint copies of the complex flag variety associated
to $G$.  In this setting, we still have a ring 
$\ring_\Gamma = \C[a_1, \dots, a_r]/I$.  However, $\ring_\Gamma$ is no
longer an integral domain, and  we cannot work with fields of fractions
analogous to $\mathcal{K}_\Gamma$; this may cause some difficulties in
generalizing the regular representation $\regr$ to this situation.  At any
rate, we do not need $\regr$ if our main goal is to prove that the
generating functions are rational.  We simply work with the trace from
$\Gamma'$ to $\Gamma$ as defined by coset representatives; this works best
if we pass to a smaller subgroup $\Gamma'' \subset \Gamma' \subset \Gamma$
such that $\Gamma''$ is a normal subgroup of $\Gamma$, similarly to taking
$\Gamma'' = \Gamma(\levelN)$ earlier in this article.  We then do our
computations in the algebra $\mathcal{A}$ obtained as a localization of 
$\ring_{\Gamma''}[A_1, \dots, A_r]$ by inverting all elements that are
congruent to $1$ modulo the ideal $(A_1, \dots, A_r)$.  (These elements are
already invertible in $\ring_{\Gamma''}[[A_1, \dots, A_r]]$, and so
$\mathcal{A}$ injects into the power series ring.)  Then the the expression
corresponding to $F = \tr^{\Gamma'}_{\Gamma} G(A_1 a_1', \dots A_r a_r')$
is a $\Gamma$-invariant expression in $\mathcal{A}$, which can be put over
a $\Gamma$-invariant common denominator (the ``norm'' of
$\prod_i (1 - A_i a'_i)$).  Then our generating function corresponding to
$F$ has a numerator and a denominator in
$\ring_{\Gamma''}[A_1, \dots, A_r]$ that are both invariant  
under $\Gamma$, since $F$ itself is invariant.  If follows that the
coefficients of the numerator and denominator of $F$ are $\Gamma$-invariant
elements of $\ring_{\Gamma''}$, i.e., elements of $\ring_\Gamma$, as
desired.





\providecommand{\bysame}{\leavevmode\hbox to3em{\hrulefill}\thinspace}
\providecommand{\MR}{\relax\ifhmode\unskip\space\fi MR }
\providecommand{\MRhref}[2]{%
  \href{http://www.ams.org/mathscinet-getitem?mr=#1}{#2}
}
\providecommand{\href}[2]{#2}

\end{document}